\documentclass[12pt]{article}

\usepackage[margin=1in]{geometry}
\usepackage{graphicx}
\usepackage{amssymb,amsmath}
\usepackage{epstopdf}
\usepackage{url}
\usepackage{amsthm}
\usepackage[normalem]{ulem}
\usepackage{dsfont}
\usepackage{amssymb}
\usepackage{graphics}
\usepackage{tikz}
\usetikzlibrary{arrows.meta}
\usepackage{array, multirow, boldline}
\usepackage{multirow}
\usepackage{xcolor}
\usepackage[hidelinks]{hyperref}
\usepackage{color,soul}
\usepackage{mathtools}
\usepackage{tablefootnote}

\usetikzlibrary{shapes.geometric,arrows,shapes,snakes,automata,backgrounds,petri}
\tikzstyle{startstop} = [rectangle, rounded corners, minimum width=1cm, minimum height=1cm,text centered, draw=black]
\tikzstyle{io} = [rectangle, rounded corners, minimum width=1cm, minimum height=1cm,text centered, draw=black]
\tikzstyle{process} = [rectangle, minimum width=3cm, minimum height=1cm, text centered, draw=black]
\tikzstyle{arrow} = [thick,->,>=stealth,black]
\tikzstyle{darrow} = [dashed,->,>=stealth, draw=blue!1000]
\tikzstyle{garrow} = [dashed,->,>=stealth, draw=violet!100]
\tikzstyle{barrow} = [dashed,->,>=stealth]
\tikzstyle{init} = [pin edge={thick, ->, black}]
\tikzstyle{int} = [pin edge={to-,thick, ->, black}]
\usepackage{wrapfig,lipsum,booktabs}

\newtheorem{prop}{Proposition}
\newtheorem{theorem}{Theorem}

\newcommand{\rv}{\mathcal{R}_v}
\newcommand{\rr}{\mathcal{R}_0}

\newcommand{\n}{\noindent}

\newcommand{\beq}{\begin{equation}}
\newcommand{\eeq}{\end{equation}}

\newcommand{\benum}{\begin{enumerate}}
\newcommand{\eenum}{\end{enumerate}}
\newcommand{\bitem}{\begin{itemize}}
\newcommand{\eitem}{\end{itemize}}

\newcommand{\ybar}{\bar{y}}

\renewcommand{\hbar}{\bar{h}}
\newcommand{\etabar}{\bar{\eta}}

\newcommand{\psibar}{\bar{\psi}}
\newcommand{\zetabar}{\bar{\zeta}}

\newcommand{\sigmabar}{\bar{\Sigma}}

\definecolor{igreen}{rgb}{0.0, 0.56, 0}

\title{Incorporating changeable attitudes toward vaccination into an SIR infectious disease model}
\author{Yi Jiang, Kristin M. Kurianski, Jane Lee, Yanping Ma, \\ Daniel Cicala, Glenn Ledder}
\date{May 2, 2024}

\begin{document}

\maketitle

\begin{abstract}
We develop a mechanistic model that classifies individuals both in terms of epidemiological status (SIR) and vaccination attitude (willing or unwilling), with the goal of discovering how disease spread is influenced by changing opinions about vaccination. Analysis of the model identifies existence and stability criteria for both disease-free and endemic disease equilibria.  The analytical results, supported by numerical simulations, show that attitude changes induced by disease prevalence can destabilize endemic disease equilibria, resulting in limit cycles.
\end{abstract}

\vspace{1in}

\section{Introduction}

Formulating public health policies toward infectious diseases requires the capability to predict how the infectious population will change over time under various scenarios.  Good modeling of the effects of public health measures is essential.  In the epidemic phase of a novel infectious disease, public health measures are limited to basic precautions such as quarantine, masking, and isolation.  These tools become less important after a disease enters the endemic phase, where a large part of the population is no longer susceptible and infectious populations are a small fraction of the total population.  From then on, the most important public health tool is vaccination, if a vaccine is available, particularly for diseases where immunity is gradually lost over time.

The traditional way to incorporate vaccination into epidemiology models is as a single-phase spontaneous transition from a susceptible class to a recovered class \cite{brauer2019mathematical,cai2018global,martcheva2015introduction}. Among the assumptions implicit in this design are that everyone will get a vaccination unless they are infected first and that there are no limits to supply and distribution.  In practice, however, a significant fraction of people do not obtain a vaccination, either because of ineligibility due to health issues or because of a refusal to receive the vaccine.  Moreover, in the early stages of a pandemic, the rate of vaccination is slowed by limited availability and distribution capacity. Ledder \cite{ledder2022puiru} addresses these issues by proposing realistic models for supply and distribution in the epidemic phase and by partitioning the susceptible class into \emph{pre-vaccinated} and \emph{unprotected} subclasses, the former comprising those who are willing and able to be vaccinated and the latter comprising those who are unwilling or unable.  In the epidemic phase, when vaccines may not yet be readily available, having an unprotected subclass makes only a slight difference because there is a pool of pre-vaccinated people who are waiting for a vaccination opportunity.  However, the existence of an unprotected subclass is an important driver of the long-term disease dynamics.

Ledder's models \cite{ledder2022puiru} assume that the fraction of unprotected individuals is fixed over time.  In reality, many people do not hold fixed attitudes toward vaccination and other public health measures, but instead react to a variety of circumstances.  Dub\'{e} et al.\ \cite{dube2013hesitancy} studied this phenomenon in light of increased hesitancy toward vaccination for the standard childhood diseases by parents in Canada.  The authors quantified individual attitudes toward vaccination on a continuum from outright refusal to full acceptance.  An individual's location within this continuum is informed primarily by their perception of the risk of the disease (and the vaccine) and the perceived importance of vaccination in mitigating that risk. These perceptions are based on the individual's background knowledge, acquired information (accurate or not), and cultural norms.  Acquired information can come from public health announcements, recommendations of health professionals, communication with acquaintances, and media posts.

Prior to the COVID-19 pandemic, there was little published work combining epidemiological dynamics with the dynamics of opinion.  An early paper on this topic examines the impact of human behavior on a disease outbreak for which control is limited to the basic mitigation procedures of quarantine and isolation \cite{funk2009spread}.  This paper does not consider these behaviors to be subject to hesitancy, but assumes that the primary driver of opinion is awareness of the disease outbreak.  The model places individuals in groups having different levels of awareness; an individual's awareness level can increase through information transmission (conversations with others) and fade over time through a natural decay process.

Other approaches to combined epidemiological and opinion dynamics have appeared since vaccine hesitancy became a significant issue in the COVID-19 pandemic.  Ali et al.\ \cite{ali2021countering} uses an epidemiological model for smallpox as part of a bioterrorism scenario.  The vaccination component categorizes people as cooperative or non-cooperative.  Opinion dynamics is driven by encounters between individuals, with no consideration of public health or social media communication, and is assumed to result in increased cooperation.  A related paper \cite{ali2023impact} couples a COVID-19 model with an opinion dynamics model similar to that of Ali et al.\ \cite{ali2021countering}, but with an additional mechanism for the influence of public communications.  The authors assume that interpersonal and social communications increase cooperation.

Another trend is to use an opinion dynamics model based on statistical physics; this formulation treats opinion as a continuous variable that changes through random interactions with other individuals and stochastic variation \cite{albi2022, bertaglia2021, dellamarca2022, dimarco2021, zanella2023}.  With no processes involving public communication or epidemiological conditions, the model outcomes are strongly dependent on the initial distribution of opinion.

The above cited papers do not include a mechanism for opinion to respond to the epidemiological state, which we take as the most quantifiable driver of opinion change.  A seminal paper by Bauch \cite{bauch2005} does include such a mechanism by assuming that the fraction of individuals who get vaccinated in a SIR model where vaccination confers immunity is subject to a dynamic game theoretic model based on the idea that individuals determine whether to accept vaccination by making a rational choice comparing the perceived benefits of vaccination and non-vaccination.
Individuals change their status when information, tracked using an \emph{information index}, is exchanged with other individuals, with an overall increase in vaccination behavior when disease prevalence exceeds a certain threshold.  d'Onofrio et al.\ \cite{donofrio2011} extend Bauch's model with more general functions for the status change dynamics, as well as functions with a time delay.  In a similar model \cite{donofrio2012}, these same authors add public information as a contributor to the information index.

Models developed by d'Onofrio et al.\ \cite{donofrio2007} and Buonomo et al.\ \cite{buonomo2013} also use an information index, but have it directly impacting the vaccination rate rather than changing rates through individuals interacting. d'Onofrio et al.\ incorporate the information index into a SIR model and the Buonomo in an SEIR model, both assuming that vaccination confers lifelong immunity.

The onset of COVID-19 has influenced the direction of mathematical epidemiology in numerous ways, including the heightened importance of individual attitudes toward vaccination and the phenomena of imperfect vaccination, in which a vaccine decreases the rate of transmission without conferring immunity.  Buonomo et al.\ \cite{buonomo2022} incorporate an information index model into a COVID-19 epidemiological model with imperfect vaccination.  The information index changes through a process of linear decay toward an equilibrium proportional to the infectious population, and the vaccination rate increases monotonically with the current value of the index.  The model is shown to have a unique endemic disease equilibrium when the disease-free equilibrium is unstable, but general stability results are not established. Another recent paper, Zuo et al.\ \cite{zuo2023} returns to the game-theoretic approach of Bauch et al.\ \cite{bauch2005} and d'Onofrio et al.\ \cite{donofrio2011}, but with a base epidemiological model that assumes imperfect vaccination.  There is no existence or stability analysis for endemic disease equilibria.

Xuan, et al.\ \cite{xuan2020network} also feature vaccination behavior dependent on disease prevalence. The authors use an agent-based SIS model with each individual's avoidance behavior based on their opinion attribute. Changes in this attribute over time are based in part on infection probability and in part on a process of consensus-building among agents connected in the social network.  Because the disease state is such an obvious driver of attitude toward vaccination, it is worthwhile to consider models in which this single driver of opinion dynamics is added to a relatively simple epidemiological model.

Our experience with COVID-19 discloses two additional features that are not included in any of the papers with vaccine hesitancy cited above.  First, vaccine hesitancy to COVID-19 has not merely slowed the rate of vaccination but has established a significant class of people who are permanently unvaccinated for a variety of reasons.  Second, vaccine hesitancy plays an enhanced role for diseases, such as COVID-19, where immunity is short-lived.

Our goal in this paper is to build a model that includes gradual loss of immunity and where vaccine hesitancy consists of an unprotected subclass in the manner of Ledder \cite{ledder2022puiru} but with transition processes that move individuals between the pro-vaccination and unprotected subclasses with rate parameters that depend on disease prevalence. To focus on the impact of these features, we make reasonable simplifying assumptions for other elements of the model, including a SIR disease structure, perfect vaccination in the sense that vaccination confers the same degree of immunity as recovery, and persistence of protection through updated vaccinations for people in the removed class who maintain a pro-vaccination attitude.

We allow the parameters that control the attitude-change processes to be nonlinear functions of disease prevalence, giving specific forms for these functions only when needed for examples.  We will see that the nonlinearity in the functions that determine the attitude-change rate parameters is sufficient to yield a high degree of instability, including limit cycles.

We begin in Section \ref{sec:modeldev} with the development of the model and its reformulation with dimensionless variables. Section \ref{sec:results} presents the analytical results for the model.  Numerical simulations in Section \ref{sec:numerical} illustrate the results with several example numerical simulations. Finally, we conclude with a discussion of the results and their practical implications in Section \ref{sec:discussion}.


\section{Model Development}\label{sec:modeldev}

We begin by presenting an initial model, which we then scale to provide a version more convenient for analysis.  The model also requires two functions to describe how the stance towards vaccination depends on the infectious population; some example functions are given at the end of this section.  Notation for the model is summarized in Tables \ref{tb:originalvars} and \ref{tb:rescaled}.

\subsection{The Initial Model}

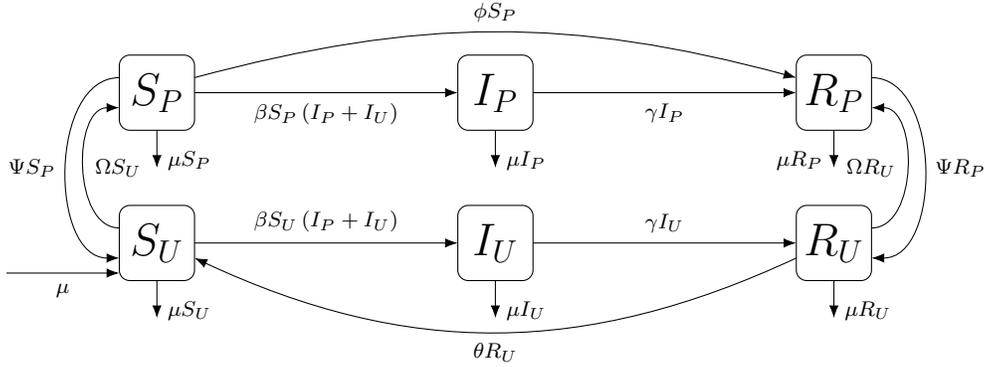
\begin{figure}
    \centering
    \[
    \begin{tikzpicture}[scale=1]
        \draw [rounded corners] (1.5-2,1.5) rectangle (2.5-2,2.5); 
        \draw [rounded corners] (4,1.5) rectangle (5,2.5); 
        \draw [rounded corners] (6.5+2,1.5) rectangle (7.5+2,2.5); 
        \draw [rounded corners] (1.5-2,-0.5) rectangle (2.5-2,0.5); 
        \draw [rounded corners] (4,-0.5) rectangle (5,0.5); 
        \draw [rounded corners] (6.5+2,-0.5) rectangle (7.5+2,0.5); 
        \draw [-Latex] (2.5-2,2) to node [below, pos=0.5] {\scriptsize $\beta S_P \left( I_P+I_U \right)$} (4,2);
        \draw [-Latex] (5,2) to node [below, pos=0.5] {\scriptsize $\gamma I_P$} (6.5+2,2);
        \draw [-Latex] (2-2,1.5) to node [right, pos=0.8] {\scriptsize $\mu S_P$} (2-2,1);
        \draw [-Latex] (4.5,1.5) to node [right, pos=0.8] {\scriptsize $\mu I_P$} (4.5,1);
        \draw [-Latex] (7+2,1.5) to node [left, pos=0.8] {\scriptsize $\mu R_P$} (7+2,1);      %
        \draw [-Latex, out=15, in=165] (2.5-2,2.2) to node [above, pos=0.5] {\scriptsize $\phi S_P$} (6.5+2,2.2);
        \draw [-Latex] (0-2,-0.4) to node [below, pos=0.5] {\scriptsize $\mu$} (1.5-2,-0.4);
        \draw [-Latex] (2.5-2,0) to node [above, pos=0.5] {\scriptsize $\beta S_U \left( I_P+I_U \right)$} (4,0);
        \draw [-Latex] (5,0) to node [above, pos=0.5] {\scriptsize $\gamma I_U$} (6.5+2,0);
        \draw [-Latex] (2-2,-0.5) to node [right, pos=0.8] {\scriptsize $\mu S_U$} (2-2,-1);
        \draw [-Latex] (4.5,-0.5) to node [right, pos=0.8] {\scriptsize $\mu I_U$} (4.5,-1);
        \draw [-Latex] (7+2,-0.5) to node [right, pos=0.8] {\scriptsize $\mu R_U$} (7+2,-1);
        \draw [-Latex, out=205, in=-25] (6.5+2,-0.2) to node [below, pos=0.5] {\scriptsize $\theta R_U$} (2.5-2,-0.2);
        \draw [-Latex, out=180, in=180] (1.5-2,0.2) to node [right, pos=0.5] {\scriptsize $\Omega S_U$} (1.5-2,1.8);
        \draw [-Latex, out=180, in=180] (1.5-2,2.2) to node [left, pos=0.5] {\scriptsize $\Psi S_P$} (1.5-2,-0.2);
        \draw [-Latex, out=0, in=0] (7.5+2,0.2) to node [left, pos=0.5] {\scriptsize $\Omega R_U$} (7.5+2,1.8);
        \draw [-Latex, out=0, in=0] (7.5+2,2.2) to node [right, pos=0.5] {\scriptsize $\Psi R_P$} (7.5+2,-0.2);
        \node () at (2-2,0) {\Large $ S_U $};
        \node () at (4.5,0) {\Large $ I_U $};
        \node () at (7+2,0) {\Large $ R_U $};
        \node () at (2-2,2) {\Large $ S_P $};
        \node () at (4.5,2) {\Large $ I_P $};
        \node () at (7+2,2) {\Large $ R_P $};
    \end{tikzpicture}
\]
    \caption{Schematic of a model incorporating pro-vaccination and unprotected subclasses with transition processes. Relevant variables and parameters are defined in Table \ref{tb:originalvars}.}
    \label{fig:schematic}
\end{figure}

\begin{table}[ht]
\caption{Definitions and dimensions of variables and parameters in the original model. }
\begin{center}
{\renewcommand{\arraystretch}{1.2}
\begin{tabular}{w{c}{1.5cm}w{c}{2cm}p{4.5in}}
\hline\hline
\textbf{Quantity} & \textbf{Dimension} & \textbf{Definition} \\
\hline
$I \;$ & dimensionless & Total Infectious population fraction \\
$I_P \;$ & dimensionless & Pro-vaccination Infectious population fraction \\
$I_U \;$ & dimensionless & Unprotected Infectious population fraction \\
$R_P \;$ & dimensionless & Pro-vaccination Recovered population fraction \\
$R_U \;$ & dimensionless & Unprotected Recovered population fraction \\
$S_P \;$ & dimensionless & Pro-vaccination Susceptible population fraction \\
$S_U \;$ & dimensionless & Unprotected Susceptible population fraction \\
$t \;$ & time & Time \\
\hline
$\beta \;$ & 1/time & Rate coefficient for infection \\
$\gamma \;$ & 1/time & Rate coefficient for recovery \\
$\theta \;$ & 1/time & Rate coefficient for immunity loss \\
$\mu \;$ & 1/time & Rate coefficient for birth and death \\
$\phi \;$ & 1/time & Rate coefficient for vaccination \\
$\Psi(I) \;$ & 1/time & Rate coefficient for Pro-vaccination to Unprotected transition \\
$\Omega(I) \;$ & 1/time & Rate coefficient for Unprotected to Pro-vaccination transition \\
\hline \hline
\end{tabular}}
\end{center}
\label{tb:originalvars}
\end{table}

Our model extends that of the PUIRU model presented by Ledder \cite{ledder2022puiru}. Our model includes additional processes that allow individuals to change their status from vaccine acceptor to non-acceptor and vice versa. Figure \ref{fig:schematic} shows a schematic diagram of our model with definitions of relevant variables and parameters in Table \ref{tb:originalvars}. The model is based on the following assumptions:

\benum
\item
The population can be divided into compartments representing the susceptible ($S$), infectious ($I$), and recovered ($R$) population fractions.  Each is divided into two subgroups: a `pro-vaccination' subclass of people who are willing and able to be vaccinated (subscript $P$) and a complementary `unprotected' subclass for those who are either unwilling or unable to be vaccinated (subscript $U$).
\item
Susceptible individuals become infectious through contact with infectious individuals with rate constant $\beta$, and infectious individuals recover through a single-phase transition process with rate constant $\gamma$.
\item
There is no disease-induced mortality. Both the birth and death rates are denoted by $\mu$ and the initial conditions are chosen so that the total population remains constant at a value of 1. Additionally, all newborns are assumed to be unprotected, and therefore, categorized within the $S_U$ subclass.
\item \label{list:assumption}
Susceptible and recovered individuals can change between pro-vaccination and unprotected subclasses, with rate coefficients $\Psi(I)$ and $\Omega(I)$, both of which are functions of the total infectious population. We neglect the possibility of status changes among the Infectious population, as people remain in that class for only a short time relative to other classes.
\item
Individuals in the pro-vaccination susceptible subclass ($S_P$) move directly to the corresponding recovered subclass ($R_P$) upon vaccination, with rate constant $\phi$.
\item
Unprotected recovered individuals ($R_U$) lose immunity over time, with rate constant $\theta$.  Pro-vaccination recovered individuals ($R_P$) remain in that subclass because they receive booster vaccination doses as needed, except for those who transfer to the other recovered subclass ($R_U$) (at rate $\Psi R_P$) because they switch from pro-vaccination to unprotected.
\eenum

We call particular attention to Assumption \ref{list:assumption}.  The rate coefficients being functions of the system's state adds additional nonlinearity, which we can surmise might produce results with features that are not seen in standard epidemic models where the only nonlinearity is in the infection process.  This assumption adds a degree of additional realism to the model.  The patterns of behavior seen in response to COVID-19 show that a lessening of the disease burden induces people to allow their immunity to lapse.

These assumptions lead to the following model,
\begin{align}
\begin{aligned}
\frac{dS_P}{dt} &= -\beta S_PI -\phi S_P-\mu S_P-\Psi S_P+\Omega S_U , \\
\frac{dS_U}{dt} &= \mu - \beta S_UI+\theta R_U - \mu S_U +\Psi S_P -\Omega S_U, \\
\frac{dI_P}{dt} &= \beta S_PI - \mu I_P - \gamma I_P, \\
\frac{dI_U}{dt} &= \beta S_UI - \mu I_U - \gamma I_U, \\
\frac{dR_P}{dt} &= \gamma I_P +\phi S_P - \mu R_P-\Psi R_P+\Omega R_U, \\
\frac{dR_U}{dt} &= \gamma I_U -\theta R_U- \mu R_U +\Psi R_P-\Omega R_U,
\end{aligned}
\end{align}
where $I=I_P+I_U$.

\subsection{The Reformulated Model}

\begin{table}
\caption{Rescaled model parameters and variables (all dimensionless and non-negative).}
\begin{center}
{\renewcommand{\arraystretch}{1.2}
\begin{tabular}{p{3cm}p{4.5in}}
\hline\hline
\textbf{Quantity} & \textbf{Definition} \\
\hline
$a$ & Exponential parameter in status-change functions $\omega_1$, $\omega_2$ \\
$d$ & Scaling parameter in status-change function $\omega_2$ \\
$h=\theta/\mu$ & Expected number of immunity losses in a lifespan \\
$I = I_P+I_U$ & Total population fraction of infectious individuals \\
$J = I_P$ & Pro-vaccination Infectious population fraction \\
$P = S_P$ & Pro-vaccination Susceptible population fraction \\
$R=R_P+R_U$ & Total population fraction of Recovered individuals \\
$S=S_P+S_U$ & Total population fraction of Susceptible individuals \\
$T=\mu t$ & Rescaled time \\
$U=S_U\!+\!I_U\!+\!R_U$ & Total population fraction of Unprotected individuals \\
$v=\phi/\mu$ & Ratio of expected lifespan to expected time for the initial vaccination\\
$X=R_U$ & Unvaccinated Recovered population fraction \\
$Y=\epsilon^{-1}I$ & Rescaled total proportion of infectious population \\
$Z=\epsilon^{-1}J$ & Rescaled Pro-vaccination Infectious population fraction \\
$\rr=\beta/(\gamma+\mu)$ & Basic reproduction number in absence of vaccination \\
$\epsilon=\mu/(\gamma+\mu)$ & Ratio of expected infectious duration to expected lifespan \\
$\omega(Y)=\Omega/\mu$ & Rescaled rate coefficient for Unprotected to Pro-vaccination transition \\
$\psi(Y)=\Psi/\mu$ & Rescaled rate coefficient for Pro-vaccination to Unprotected transition  \\
\hline \hline
\end{tabular}}
\end{center}
\label{tb:rescaled}
\end{table}

We rescale the model to simplify the stability analysis. The dependent variables are already scaled as fractions of the total population.  We introduce a dimensionless time $T$ using the population time scale $1/\mu$, that is
$T=\mu t$; hence,
\[
\frac{d}{dt} = \mu \frac{d}{dT}.
\]
We also introduce dimensionless parameters
\begin{equation}\label{eq:dimlessparams}
\epsilon = \frac{\mu}{\gamma+\mu}, \quad \mathcal{R}_0=\frac{\beta}{\gamma+\mu}, \quad h=\frac{\theta}{\mu}, \quad v = \frac{\phi}{\mu}, \quad \omega = \frac{\Omega}{\mu}, \quad \psi = \frac{\Psi}{\mu}
\end{equation}
and use $'$ to denote the derivative with respect to $T$.  The intuitive meanings of all newly-defined notations here are included in Table~\ref{tb:rescaled}.

The system is a little simpler if the combined variables $S$, $I$, and $R$ are used instead of the subclasses $S_U$, $I_U$, and $R_P$.  We also simplify the notation by using non-subscripted variables $P$, $J$, and $X$ in place of $S_P$, $I_P$, and $R_U$, respectively.  Finally, the equations are reordered to facilitate the use of the next generation matrix method for the calculation of the basic reproduction number in the presence of vaccination, $\rv$.  With all of these changes, the model becomes

\begin{align}
\begin{aligned}
\label{syst1}
I' &= \epsilon^{-1} (\rr S-1)I, \\
J' &= \epsilon^{-1} (\rr PI-J ),\\
S' &= 1 - S + h X -v P -\epsilon^{-1} \rr SI, \\
P' &= -(1+v) P -\psi(Y) P +\omega(Y) (S-P)- \epsilon^{-1}\rr PI, \\
R' &= \left( \epsilon^{-1}-1 \right) I + v P - h X - R, \\
X' &= \left( \epsilon^{-1}-1 \right) (I-J) - (1+h) X - \omega(Y)
X+ \psi(Y) (R-X).
\end{aligned}
\end{align}

The system \eqref{syst1} incorporates the naive assumption that all dependent variables are $O(1)$ quantities on the time scale of basic population dynamics.  The resulting system is not correctly scaled for analysis of long-term behavior,
as the $R$ equation would then imply $I=O(\epsilon)$ at equilibrium, contradicting the inherent assumption that all variables are $O(1)$.  Other equations are similarly inconsistent.  The solution is to build in the assumption that infectious populations are in fact $O(\epsilon)$ on this time scale \cite{ledder2023using}. Note that $\epsilon$ represents the ratio of expected time spent in the infectious compartment to mean lifetime. Assuming infections last on the order of weeks and a lifetime is roughly 80 years, this parameter is indeed small. We therefore replace $I$ and $J$ by
\beq
I = \epsilon Y, \qquad J=\epsilon Z.
\eeq
At the same time, we introduce some lumped quantities that are convenient now or will become so later:
\benum
\item
As needed, a bar over the top of a quantity indicates that quantity plus 1; for example, $\sigmabar=\Sigma+1$.
\item
It will be convenient to remove a factor of $\rr^{-1}$ from some of the equilibrium values; in anticipation of this, we define
\beq
\label{eqvars}
y=\rr Y, \qquad z=\rr Z, \qquad p=\rr P, \qquad x=\rr X.
\eeq
\item
Various sums of parameters arise in the analysis:
\beq
\label{params1}
\Sigma=\psi+\omega, \qquad \zeta=\sigmabar+v, \qquad \eta=\sigmabar+h.
\eeq
\eenum

\begin{table}
\caption{Modified and grouped variables and parameters used in analysis.}
\begin{center}
{\renewcommand{\arraystretch}{1.2}
\begin{tabular}{p{4cm}p{4.5in}}
\hline\hline
\textbf{Quantity} & \textbf{Definition} \\
\hline
$\Sigma=\psi+\omega $ & Expected number of opinion changes in lifespan \\
\hline
$\bar{h}=h+1$ & Parameters with bars indicate addition by $1$\\
$y=\rr Y, \cdots$ & Lower case model variables correspond to multiplication by $\rr $; Applies to $Y$, $Z$, $P$, and $X$, but not $S$, $R$, or $U$\\
\hline
$\rv$ & Basic reproduction number in the presence of vaccination \\
$r=\rr-1$ & Parameter grouping defined in Proposition \ref{DFE} \\
$\zeta = \bar{\Sigma}+v$ & Parameter grouping defined in \eqref{params1} \\
$\eta = \bar{\Sigma}+h$ & Parameter grouping defined in \eqref{params1} \\
$\xi = \bar{\Sigma}+h\bar{\psi}$ & Parameter grouping defined in Proposition \ref{DFE} \\
$\rho =\sigmabar(v+rh)+\rr hv$ &Parameter grouping defined in Proposition \ref{DFE} \\
$\chi = \bar{\Sigma}+v+y$ & Parameter grouping defined in Proposition \ref{EDEomegastar} \\
\hline \hline
\end{tabular}}
\end{center}
\label{tb:grouped}
\end{table}

With these changes, we obtain the final form of the system (note that the lumped parameters $\eta$ and $\zeta$ are functions of $Y$ because each contain a $\Sigma$ term):

\begin{align}
\begin{aligned}
\label{syst2}
\epsilon Y' &= (\rr S-1)Y, \\
\epsilon Z' &= \rr PY-Z, \\
S' &= 1-S-vP+hX-\rr SY, \\
P' &= \omega(Y) S-\zeta(Y) P-\rr PY, \\
R' &= (1-\epsilon)Y+vP-R-hX, \\
X' &= (1-\epsilon)(Y-Z)+\psi(Y)R-\eta(Y) X.
\end{aligned}
\end{align}

\subsection{The Status-Change Functions $\omega$ and $\psi$}

Regardless of whether the unprotected-to-protected transition rate coefficient $\omega$ is increasing or decreasing, it seems reasonable that the protected-to-unprotected transition rate coefficient $\psi$ should be moving in the opposite direction.  While one of these functions might be more sensitive to the infectious population fraction, there is no way to know which. As a first approximation, it is not unreasonable to think that the sum $\Sigma=\omega+\psi$ might be approximately constant.  To keep the model as simple as possible for the elucidation of general properties, we will think of $\Sigma$ as constant whenever we are ready to consider specific examples, specifying $\omega$ and replacing $\psi$ by $\Sigma-\omega$.

We can identify two different classes of plausible status-change functions $\omega$, given the assumption of constant $\Sigma$.  When infection rates are low, $\omega$ should be increasing with $Y$.  This is something that typically occurs with diseases; a lower incidence leads to less strenuous vaccination campaigns, as is currently done with some vaccinations administered only to people traveling to a country with known cases.  As the infection rate increases, two types of behavior seem possible.  Healthcare authorities will undoubtedly increase the emphasis they put on the need for vaccination, and an informed public that is trusting of healthcare authorities will be more likely to be vaccinated.  On the other hand, public skepticism could lead to a situation where some people interpret higher disease burdens as vaccine failures and become less likely to be vaccinated.  Thus, we consider functions that are monotone-increasing and functions that increase to a local maximum before declining.

For elucidation of general behavior, the specific families of functions used for the monotone and non-monotone cases can be chosen for mathematical convenience.  In the sequel, we will use
\beq
\label{omega1}
\omega_1(Y;\Sigma,a) = \Sigma \left( 1- e^{-aY} \right), \qquad a>0,
\eeq
as a family of monotonic functions and
\beq
\label{omega2}
\omega_2(Y;\Sigma,a,d) = \Sigma daYe^{1-aY}, \qquad a,d>0,
\eeq
as a family of non-monotonic functions.  The latter family has the simple property that it vanishes to 0 as $Y \to \infty$; this seems reasonable, as having nearly everyone infected would likely make everyone lose confidence in the vaccine.

\subsection{Choice of Parameter Values}
\label{sect:paramvalues}

Values for parameters are problematic when a model is intended to be general rather than to match a specific disease.  It is best to defer choosing values until absolutely necessary.  In this case, we can obtain analytical results using only the approximation $\epsilon \to 0$, but parameter values are still required for visualization of results.  When this is required, we choose fixed values of $\rr=4$ and $v=50$.  The basic reproduction number does not vary greatly among diseases.  The value 4 is intermediate between diseases such as COVID-19, which are highly transmissible,\footnote{The original strain had $\rr \approx 5.7$ \cite{sanche2020high}, and more recent variants are surely higher.} and less infectious diseases such as the flu, with $\rr \approx 1.6$ \cite{biggerstaff2014estimates}.  The value $v=50$ corresponds to the assumption that the time required for a person willing to be vaccinated to obtain vaccination is about 2\% of a normal lifespan, which is about 20 months.  We should expect a higher value in practice, but the results are not sensitive to further increases in this parameter.

The parameter $h$ represents the expected number of times immunity is lost during a normal lifespan.  Actual values of this parameter can vary from 0 for diseases such as measles, where immunity tends to be permanent, to diseases such as the flu, where immunity wanes quickly.  We will generally use $h=0$ and $h=10$ as representative of diseases with long and short immunity durations, respectively.  While the parameter $\Sigma=\omega+\psi$ is variable in the general case, it is reasonable to think it might in practice be a fixed parameter that measures the flexibility of opinion in the community.  We will use $\Sigma=20$ as a practical upper bound for examples; this value represents a population in which opinions change every four years on average.  Of course a real population will consist of some people whose opinions are fixed for life and others whose opinions change quickly.

\section{Analysis and Results}\label{sec:results}

We begin with the disease-free equilibrium (DFE), continue with general results for the endemic disease equilibrium (EDE), and then examine results for some specific ideology-change functions $\psi$ and $\omega$.  Because this model is intended to study general features of epidemics, rather than reproducing data or making predictions for a specific disease, our primary interest is in qualitative results rather than quantitative details.  For this reason, it makes sense to make use of $\epsilon \to 0$ asymptotics for the endemic disease equilibrium analysis.\footnote{Analysis of the disease-free equilibrium will not require asymptotic approximation.}  Hence, terms in a sum that are $O(\epsilon)$ compared to the leading order term will be discarded.  In addition to the $O(\epsilon)$ terms in the differential equations, there will be $O(\epsilon)$ terms that arise in the stability analysis.

\subsection{The Disease-Free Equilibrium and Vaccine-Reduced Reproduction Number}

Computation of the vaccine-reduced reproduction number and identification of the stability requirement for the disease-free equilibrium are accomplished by standard methods.  We present the result here and give details in Appendix \ref{appendix:DFEprop}.

\begin{prop}
\label{DFE}
The vaccine-reduced reproduction number for the model \eqref{syst2} is
\beq
\rv =
\frac{\zeta \xi}{\zeta \xi + v \omega \eta} \, \rr,
\eeq
where $\rr$ is given by \eqref{eq:dimlessparams}, $\zeta$ and $\eta$ are given by \eqref{params1}, and $\xi=\sigmabar+h\psibar$.  The disease-free equilibrium is stable if and only if $\rv<1$; equivalently, there is a minimum value of $\omega$ needed for stability:
\beq
\label{omegacr}
\omega(0) > \omega_{cr} \coloneqq \frac{ r\hbar \sigmabar(0) \zeta(0)}{\rho(0)},
\eeq
where $\rho(Y)=\sigmabar (Y)(v+rh)+\rr hv$ and $r=\rr-1$
\end{prop}

\begin{figure}[ht]
    \centering
    \includegraphics[width=0.5\textwidth]{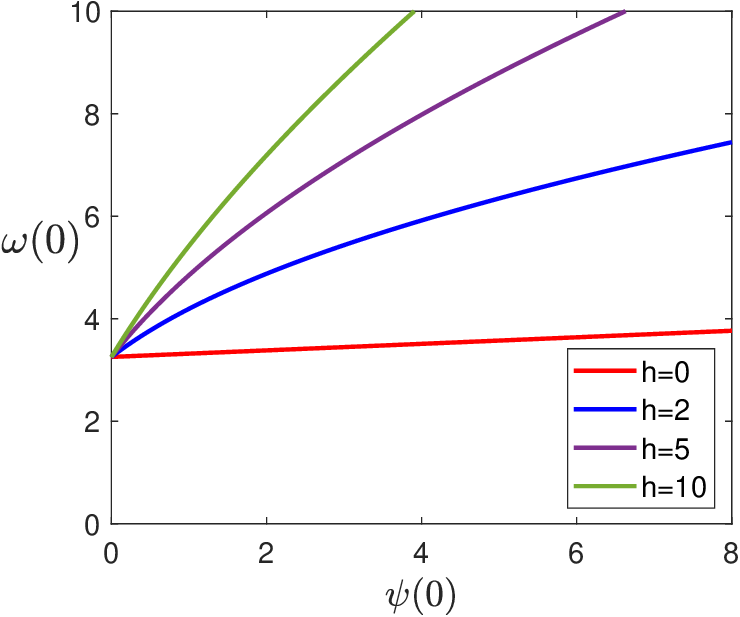}
    \caption{Stability boundaries for the disease-free equilibrium (DFE), with $h=0, 2, 5, 10$ (bottom to top) and $\rr=4$, $v=50$; the DFE is stable if the point in the $\psi \omega$-plane is above the boundary curve for given $h$.}
\label{fig:DFEstability}
\end{figure}

Figure \ref{fig:DFEstability} illustrates the stability requirement for the DFE.  Larger $\psi$ requires larger $\omega$ to achieve stability, and larger $h$ increases the $\omega$ requirement.
As a practical matter, we should expect that this stability requirement cannot be achieved, as this would require a significant effort to vaccinate young children in spite of having no active disease in the population.

\subsection{The Endemic-Disease Equilibrium}

The parameters $\omega$ and $\psi$, and the quantities derived from them, are functions of $Y$.  Thus, when we solve the equilibrium system for $y$, the implicit result that we get will be an equation that contains unknown values for $\omega$ and $\psi$ as well as $y$.  This equation can be solved for $y$ only after the functions $\omega(Y)$ and $\psi(Y)$ have been defined.  Given this situation, and the fact that the combination $\Sigma=\omega+\psi$ occurs frequently, it will be most convenient to think of the equilibrium relation as specifying a ``target'' value of $\omega$ that must be achieved at a prescribed value of $y>0$ for that value to be an equilibrium, given known values of $\Sigma$, $\rr$, $v$, and $h$.  We summarize the results here and defer details to Appendix \ref{appendix:omegaprop}.

\begin{prop}
\label{EDEomegastar} $ $
\benum
\item
For any given values of the parameters $\Sigma$, $\rr$, $v$, and $h$, endemic disease equilibria occur when
\beq
\omega(Y)=\omega^*(\rr Y),
\eeq
where
\beq
\omega^*(y) = \chi(y)p(y), \qquad p(y)=\frac{\sigmabar(r \hbar -y)}{\rho+\rr hy}, \qquad \chi(y)=\zeta+y,
\eeq
with $r=\rr-1$, $\hbar=h+1$, $\zeta=\sigmabar+v$, $\rho=\sigmabar(v+rh)+\rr hv$, and $y=\rr Y$.

\item
The function $\omega^*$ has the properties
\[ \omega^*(0)=\omega_{cr}, \qquad \omega^*(r\hbar)=0, \qquad {\omega^*}'(y)=-\sigmabar \;\frac{\kappa+2\rho y+\rr hy^2}{(\rho+\rr hy)^2}, \]
where $\kappa=\zeta \rho -\sigmabar r \hbar(v-h)$.
\item
The
condition
\beq
\label{monotone}
(1-h \zeta)\rr<\sigmabar \zeta+1
\eeq
is sufficient to guarantee that $\omega^*$ is monotone decreasing.
\eenum
\end{prop}

\n The condition \eqref{monotone} is not very restrictive.  Either of the simpler conditions
\[ \rr<v+2, \qquad h>\frac{1}{v+1} \]
is sufficient, although not necessary, to guarantee  \eqref{monotone}.  Any realistic disease scenario with available vaccine should have $v>\rr$.

\begin{figure}[ht]
    \centering
    \includegraphics[width=0.5\textwidth]{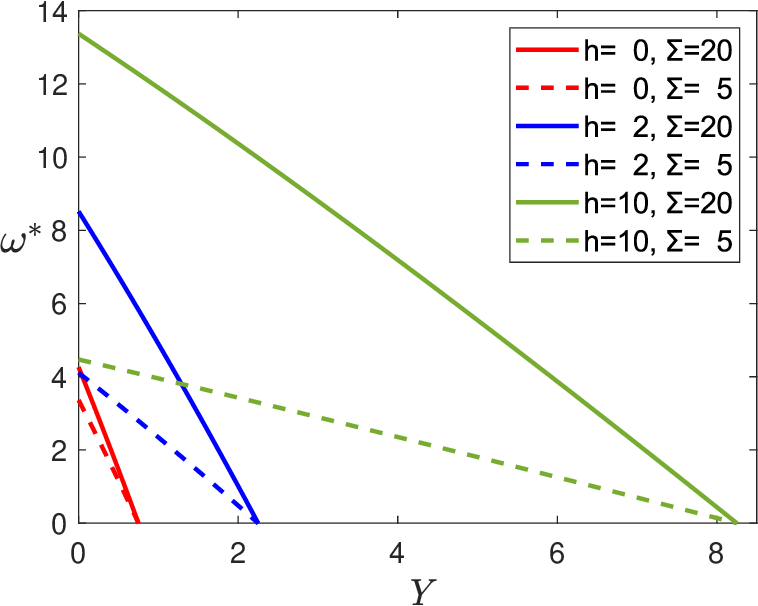}
    \caption{Examples of $\omega^*(\rr Y)$ with $\rr=4$ and values of $\Sigma$ and $h$ as appear in the legend.}
    \label{fig:omegastar}
\end{figure}

Figure \ref{fig:omegastar} shows some examples of $\omega^*(\rr Y)$.  These curves represent the values that would have to be achieved by $\omega(Y^*)$ in order for $Y^*$ to be an endemic disease equilibrium.  Monotone increasing functions for $\omega$ will have a unique EDE if and only if the disease-free equilibrium is unstable, that is, $\omega(0)<\omega_{cr}$. Endemic disease equilibria must have $y<r\hbar$, or
\beq
\label{ymax}
Y^* < Y_{max} = \left( 1-\rr^{-1} \right) \hbar.
\eeq

\subsubsection{Endemic Disease Equilibrium Stability}

While it is not normally possible to obtain analytical stability criteria for a six-component system with unspecified parameter values, we are able to do so with our system in the limit $\epsilon \to 0$ (see \cite{ledder2023using} for a description of the method).  We summarize the results here, leaving details for Appendix \ref{appendix:EDEstab}.

\begin{prop}
\label{EDEstability}
Let $Y>0$ be an equilibrium point for the system.  Given the values of $\omega$, $\psi$, and their derivatives at the equilibrium point, the stability requirements for the endemic disease equilibrium in the asymptotic limit $\epsilon \to 0$ reduce to
\beq
\label{stab1}
v\omega'-hr\psi'+(hx-vp)\Sigma' < \rr w_1,
\eeq
\beq
\label{stab2}
- \left[ v\omega'-hr\psi'+(hx-vp)\Sigma' \right] < \rr \left[ \sigmabar+w_2 \right],
\eeq
\beq
\label{stab3}
-\left[ (v \eta+hy) \omega'-rh\chi \psi'+(hx\chi -v\eta p-hyp)\Sigma' \right] < \rr \sigmabar w_2 ,
\eeq
where
\[ w_1 = h(1-p)+vp+\ybar, \qquad w_2 = v(1-p)+hp+y+\sigmabar. \]
For the simplified case where $\Sigma$ is independent of $Y$, these conditions reduce to
\beq
\label{stab}
\rr {\omega^*}'(\rr Y) < \omega'(Y) < \rr \Upsilon(\rr Y),
\eeq
where
\beq
\label{upsilon}
\Upsilon(y)=\frac{w_1(y)}{v+rh}.
\eeq
\end{prop}

\begin{figure}[ht!]
    \centering
    \includegraphics[width=0.32\textwidth]{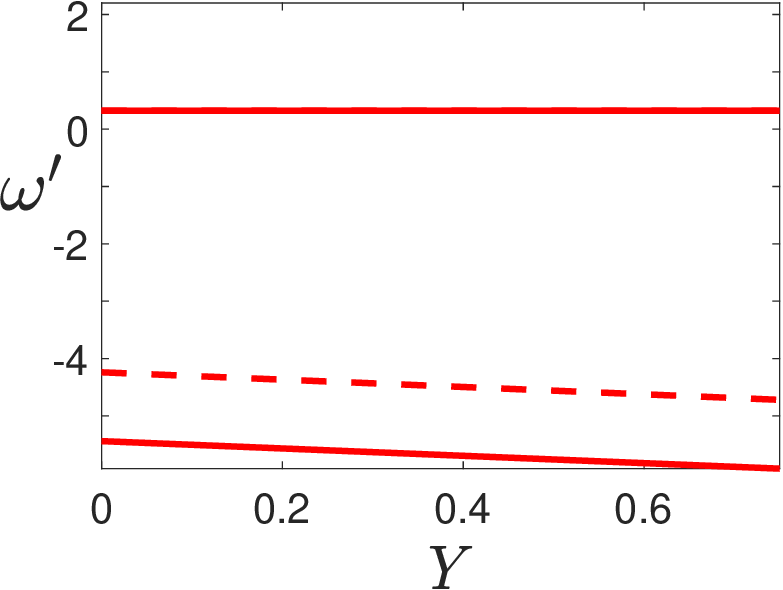}
    \hfill
    \includegraphics[width=0.32\textwidth]{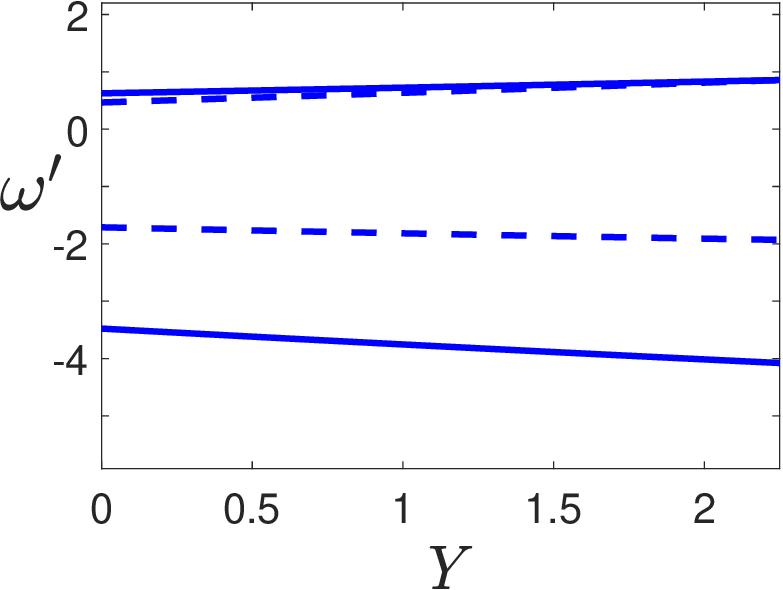}
    \hfill
    \includegraphics[width=0.32\textwidth]{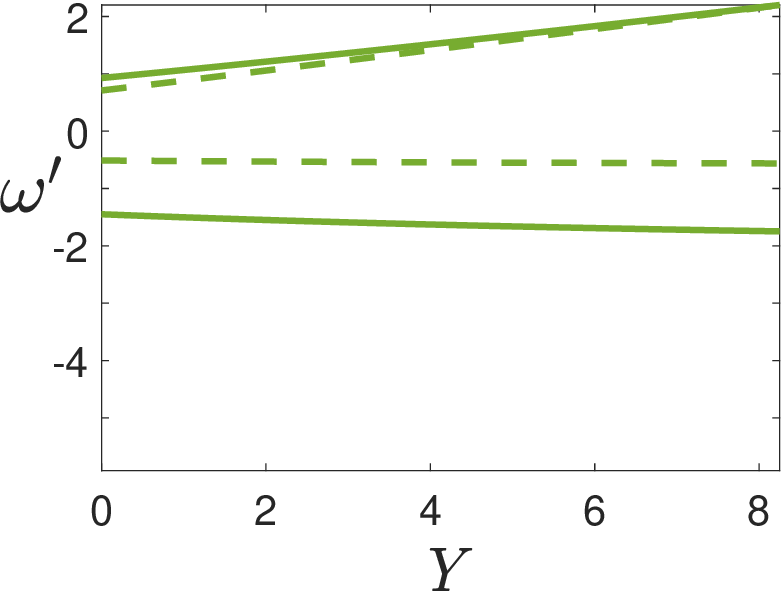}
    \caption{Stability boundaries with $\rr=4$, $v=50$; $h=0, 2, 10$ (left to right); solid: $\Sigma=20$, dashed: $\Sigma=5$. Values of $\omega'$ between the curves correspond to stable EDE condition \eqref{stab}.}
    \label{fig:stability}
\end{figure}

Figure \ref{fig:stability} illustrates the stability requirements for $\omega'$ in terms of $Y$, with $\rr=4$, $v=50$, $h=0, 2, 10$ (left to right) and $\Sigma=20$ (solid) and $\Sigma=5$ (dashed).  The constant $\omega$ case corresponds to a model with no dependence of rates on disease level, for which the (unique by Proposition \ref{EDEomegastar}) EDE will always be stable.  Instability results from greater sensitivity of $\omega$ to $Y$, with stricter stability boundaries for smaller $\Sigma$.  Interestingly, increasing $h$ makes the stability boundary more strict when $\omega'$ is negative and less strict when it is positive.

\subsection{Results for Some Specific Status-Change Functions}

We now consider results for some specific status change functions, both assuming $\Sigma$ is constant.

\subsubsection{Monotone Functions}
\begin{figure}[h!]
    \centering
    \includegraphics[width=1\textwidth]{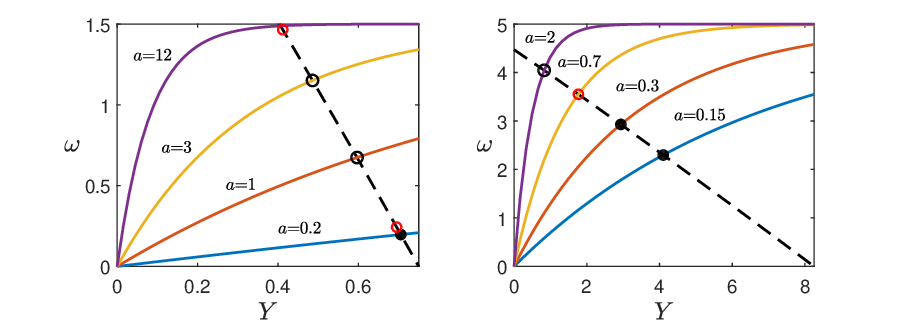}
    \caption{Examples of functions $\omega=\Sigma \left( 1-e^{-aY} \right)$ for different values of $a$ (solid) along with $\omega^*$ (dashed); parameters are $\rr=4$, $v=50$, with $\Sigma=1.5$, $h=0$ (left) and $\Sigma=5$, $h=10$ (right). Stable equilibria are marked with solid disks and unstable with an open circle, and the stability boundaries are red circles.}
\label{fig:Omegaplots1}
\end{figure}

We begin with the broad class of monotone increasing functions for $\omega(Y)$.  If $\omega(0)>\omega_{cr}$, then the equilibrium condition $\omega(Y)=\omega^*(\rr Y)$ will never be achieved and the disease-free equilibrium will be stable.  Instead, we assume $\omega(0)<\omega_{cr}$.  In this case, the combination of properties $\omega'>0$, $\omega(0)<\omega_{cr}$, and $\omega^*(r\hbar)=0$ guarantee that there will be a unique endemic disease equilibrium $Y>0$, and that it will be stable if and only if
\beq
\omega'(Y) < \rr \Upsilon(\rr Y), \qquad \Upsilon(y)=\frac{h(1-p)+vp+\ybar}{v+rh}.
\eeq

As a simple example of this class, we consider the one-parameter family given in \eqref{omega1}.  Several functions in this family are shown in Figure \ref{fig:Omegaplots1}, along with the curve $\omega^*(\rr Y)$ with $\rr=4$, $v=50$, and either $h=0$, $\Sigma=1.5$ or $h=10$, $\Sigma=5$.  The unique equilibrium for each curve is marked according to its stability result.  For $h=10$, the situation is relatively straightforward.  The critical case is $a=0.7$; larger values yield a larger slope at the intersection and mark an unstable equilibrium, while smaller values yield a stable equilibrium.  For $h=0$, the situation is a little more complicated.  Here the slope is small enough for stability with either a very small or a very large $a$; intermediate values of $a$ have a steeper slope at the intersection and yield instability.

\begin{figure}[h!]
    \centering
    \includegraphics[width=1\textwidth]{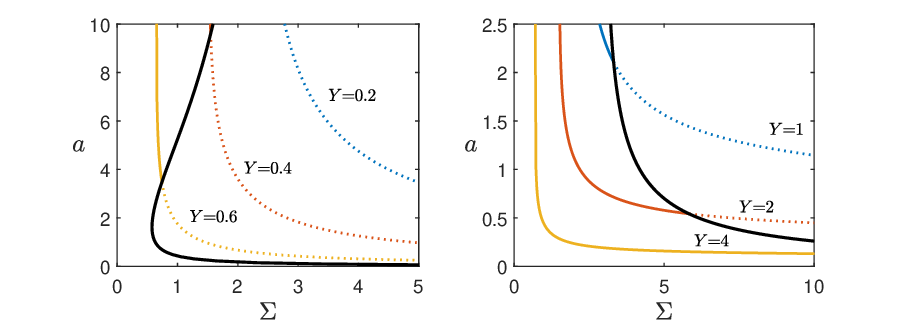}
    \caption{Examples of level curves of $Y$ in the $\Sigma a$-plane for $\omega=\Sigma \left( 1-e^{-aY} \right)$, along with the stability boundary (black), with $\rr=4$, $v=50$; $h=0$ (left), $h=10$ (right). The level curves are solid when the equilibrium is stable and dotted when the equilibrium is unstable.}
    \label{fig:EDEplots1}
\end{figure}

A broader summary of the $h=0$ and $h=10$ cases, with $\rr=4$ and $v=50$, appears in Figure \ref{fig:EDEplots1}, which shows level curves of $Y$ in the $\Sigma a$-plane along with the stability boundary.  With $h=10$, the stability boundary curve always moves to the left as $a$ increases; thus, larger $a$ is always destabilizing.  In contrast, with $h=0$ the curve moves first to the left and then to the right; here there are $\Sigma$ values for which instability occurs in an intermediate band of $a$ values.

\subsubsection{Functions with a Local Maximum}

We now turn to the broad class of functions for $\omega(Y)$ that increase to a maximum at a point $Y_m<Y_{max}$, and then decrease as $Y$ increases further.

\begin{figure}[ht]
    \centering
    \includegraphics[width=0.5\textwidth]{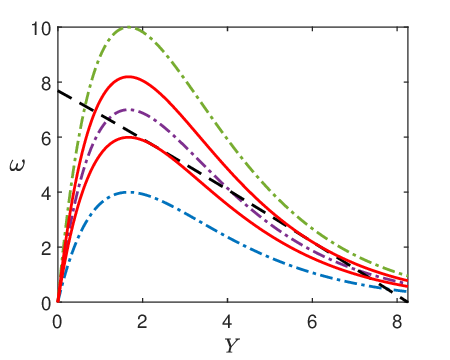}
    \caption{$\omega(Y)=\Sigma daY e^{1-aY}$, along with $\omega^*(\rr Y)$ (dashed), with $\Sigma=10$, $a=0.6$, $h=10$, $\rr=4$, and $v=50$; the red solid curves have $d=0.60$ (lower) and $d=0.82$, while the dash-dot curves are for $d=0.4, 0.7, 1.0$, bottom to top.}
    \label{fig:Omegaplots2}
\end{figure}

For a specific example, we consider the family \eqref{omega2}, which increases to a maximum value of $\Sigma d$ occurring at $aY=1$.  The requirement that  $\psi$ be non-negative translates into a parameter restriction $d \le 1$.  A sampling of this family is illustrated in Figure \ref{fig:Omegaplots2}.  The dashed curve is the equilibrium condition $\omega(Y)=\omega^*(\rr Y)$, with $\Sigma=10$, $h=10$, $\rr=4$, and $v=50$.  All curves for $\omega$ have $a=0.6$. The three dash-dot curves have $d$ values of 0.4, 0.7, and 1.0, from bottom to top.  The first of these has one equilibrium, with a value of $Y$ that is close to $Y_{max}$, where $\omega^* (\rr Y)$ is 0.  The equilibrium occurs at a point where $\omega$ is decreasing.  The last of the solid curves also has one equilibrium, this time with a small value of $Y$ and with $\omega'>0$ at the equilibrium value. The intermediate curve has three equilibria.

In general, suppose $\omega(Y)$ increases from a point $\omega(0)<\omega_{cr}$ to a point $\omega(Y_m)$? and then decreases, with $\omega(Y^*)>0$. Then the equilibrium condition $\omega(Y)=\omega^*(\rr Y)$ will necessarily be achieved at least once, but may be achieved multiple times, with bifurcation points separating curves in the family that have different numbers of equilibria.  The $\omega$ functions that mark bifurcations satisfy the requirement that the $\omega(Y)$ and $\omega^*(\rr Y)$ curves are tangent at a point, that is
\beq
\label{bifeqns}
\omega(Y)=\omega^*(\rr Y), \qquad \omega'(Y)=\rr {\omega^*}'(\rr Y).
\eeq
Because ${\omega^*}'<0$, bifurcation points can only occur in the domain where $\omega$ is decreasing, that is, when $Y>Y_{m}$; hence, each bifurcation point in the $ad$-plane corresponds to a unique $Y$ value in the interval $(Y_m, Y_{max})$, where $Y_{max}$ is defined in \eqref{ymax}.
The function $\omega_2$ in \eqref{omega2} has two bifurcation points for any given value of $a$. These are the two dotted curves in Figure \ref{fig:Omegaplots2}, corresponding to $d=0.60$ and $d=0.82$.

Comparison with earlier results discloses an interesting synergy, which holds for any function with the properties given above:
\begin{prop}
\label{multistability}
Suppose $\omega(Y)$ increases to a maximum at a value $Y_m$, and then decreases, with $\omega(0)<\omega_{cr}$ and $\omega(Y_{max})>0$.  Then for the region of solutions $Y>Y_m$, the stability boundary \eqref{stab} coincides with the bifurcation curve \eqref{bifeqns}.
\end{prop}

\n The proposition follows immediately from the fact that $\omega'<0<\rr \Upsilon(\rr Y)$ when $Y>Y_m$, thus automatically satisfying the upper bound criterion for stability.

\begin{figure}[h!]
    \centering
    \includegraphics[width=1\textwidth]{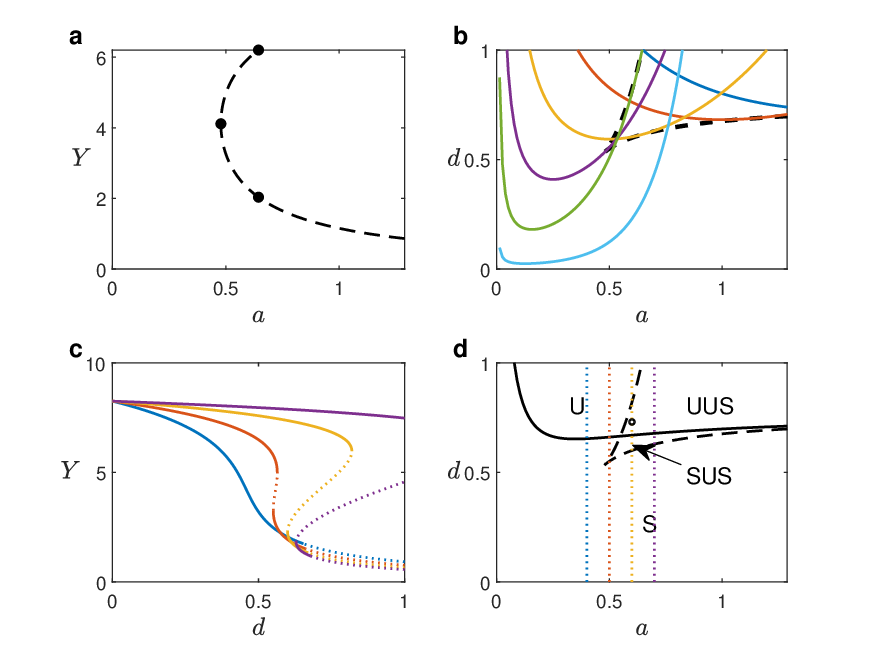}
    \caption{Properties of the system with $\omega(Y)=\Sigma da Y e^{1-aY}$, with $\Sigma=10$, $h=10$, $\rr=4$, $v=50$; \textbf{a}: The curve of bifurcation points; \textbf{b}: bifurcation curves (dashed) and level curves for $Y=0.6, 1, 2, 4, 6.4, 8$ (top to bottom in the region to the left of the bifurcation curve); \textbf{c}: solution curves for $a = 0.4, 0.5, 0.6, 0.7$ (dotted where unstable); \textbf{d}: bifurcation curve (dashed) and stability boundary (solid); regions are labeled as stable or unstable where there is a unique solution, and UUS where only the largest of three solutions is stable; the small unlabeled region in the center right is SUS; dashed lines mark the $a$ values used in \textbf{c}. The small circle marks the point used for Figure \ref{fig:nonmonotone-examps}.
    \label{fig:Bifplots2}}
\end{figure}

Figure \ref{fig:Bifplots2} illustrates the behavior of the system.  The bifurcation points are indicated by the curve in panel \textbf{a} and the dashed curves in panel \textbf{b}.  The former shows that each $Y<Y_{max}$ corresponds to a unique point in the $ad$-plane, although only those with $d<1$ meet the model requirements.  Panel \textbf{b} also shows level curves for some solution values; the curves are consistent with the bifurcation curves, as there are multiple solutions only in the region in the upper right corner of the $ad$-plane.  Panel \textbf{c} shows solution curves in the $dY$-plane, with each curve corresponding to a value of $a$ in a different range (marked with dotted curves in panel \textbf{d}).  Panel \textbf{d} also shows the bifurcation curves, along with the (black) stability boundary curve for positive $\omega'$, given by
\[ \omega'(Y) = \rr \Upsilon(\rr Y). \]
For points above this stability curve, the smallest $Y$ value (the unique $Y$ outside of the region of multiple solutions) is unstable.  When there are multiple solutions, the one with largest $Y$ is stable and the one with intermediate $Y$ is unstable, from Proposition \ref{multistability}.  The regions are labeled with their stability pattern from low $Y$ to high, with the slender unlabeled region in the middle as SUS, that is, three equilibria that are stable, unstable, and stable, from smallest to largest.  The $a$ values used in panel \textbf{c} are marked by the dotted vertical lines in panel \textbf{d}.  The first is in the region where the solution is always unique.  The second is in a very narrow region where increasing $d$ changes the stability pattern from S to SUS, back to S, and then to UUS.  The third region has the sequence S, SUS, UUS, U, and the fourth region has S, SUS, UUS.  The four solution curves in panel \textbf{c} correspond to these results.

\section{Simulations}
\label{sec:numerical}

\begin{figure}[h!]
    \centering    \includegraphics[width=1\textwidth]{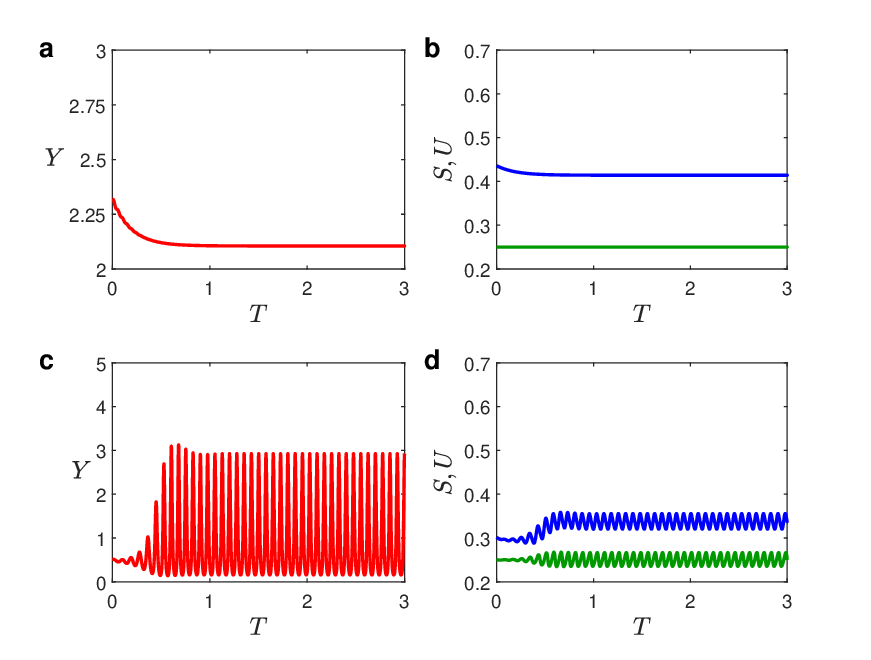}
    \caption{Examples of stable EDE (top row) and unstable EDE (bottom row) solutions for $\rr=4$, $v=50$, $h=10$, $\Sigma=2$ (top) or $\Sigma=5$ (bottom) $,\epsilon=5*10^{-4}$, and $\omega_1$ defined by \eqref{omega1} with $a=1$. The plots show the rescaled infectious fraction $Y$, the susceptible population fraction $S$ (green), and the total unprotected population fraction $U$ (blue).  Initial conditions were chosen to be near the equilibrium values.
    \label{fig:monotone-examps}}
\end{figure}

Some simulations for the monotone example function $\omega(Y)= \Sigma \left( 1- e^{-aY} \right)$ are shown in Figure \ref{fig:monotone-examps}.  The simulations use our standard values for most parameters, as explained in Section \ref{sect:paramvalues}, but using two different values for the population flexibility parameter $\Sigma$.  As predicted in Figure \ref{fig:EDEplots1}, the smaller value of $\Sigma$ results in a stable EDE with a $Y$ value near 2, while the larger results in an unstable EDE.  We see from the simulation that the latter case produces a stable limit cycle.  Note that the amplitudes of the oscillations for $S$ and $U$ are small, whereas the amplitude of the oscillations for $Y$ is large, corresponding to the classification of $Y$ as a fast variable (rate of change $O(\epsilon^{-1})$).  When $Y$ hits its peak, the susceptible population decreases and the increase in the pro-vaccination attitude causes the unprotected population to decrease as well.  These changes decrease the overall population susceptibility, which reverses the rapid increase in prevalence.

\begin{figure}[h!]
    \centering    \includegraphics[width=1\textwidth]{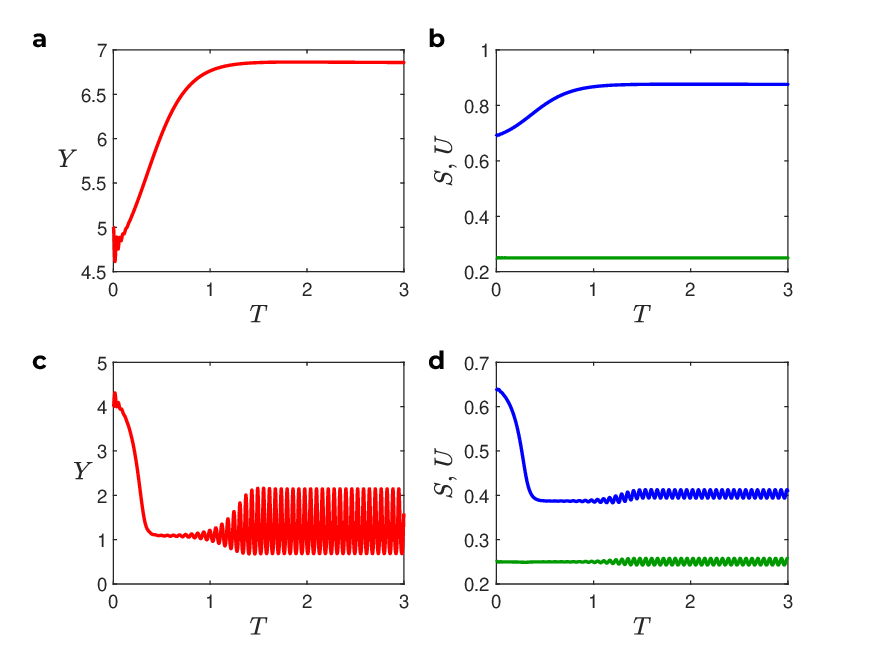}
    \caption{Examples of stable EDE (top row) and unstable EDE (bottom row) solutions for $\rr=4$, $v=50$, $h=10$, $\Sigma=10$, $\epsilon=5*10^{-4}$, and $\omega_2$ defined by \eqref{omega2} with $a=0.6$ and $d=0.73$. The plots show the rescaled infectious fraction $Y$, the susceptible population fraction $S$ (green), and the total unprotected population fraction $U$ (blue).  \textbf{Top row:} Initial conditions  $(Y_0,\,Z_0,\,S_0,\,P_0,\,R_0,\,X_0) = (5,\,0.05,\,0.25,\,0.01,\,0.75,\,0.44)$. \textbf{Bottom row:} Initial conditions  $(Y_0,\,Z_0,\,S_0,\,P_0,\,R_0,\,X_0) = (4,\,0.05,\,0.25,\,0.01,\,0.75,\,0.4)$.
    \label{fig:nonmonotone-examps}}
\end{figure}

For the non-monotone case, we focus on a single simulation experiment designed to verify the features predicted in Figure \ref{fig:Bifplots2}\textbf{d} for the parameter set corresponding to the marked circle.  The prediction is for a set of three equilibria with relatively small, medium, and large $Y$ values, observable in panel \textbf{c} as the points on the gold curve corresponding to $d=0.73$ and having the relatively uncommon result in which the largest of the three is stable and the other two unstable.  The two simulations of Figure \ref{fig:nonmonotone-examps} use initial conditions that are near the unstable equilibrium with medium $Y$ value.  The top two panels show the results of a simulation that converges to the stable large $Y$ equilibrium, while the bottom two panels show a simulation that converges to a stable limit cycle corresponding to the small $Y$ equilibrium.  If we were to keep the same value of $a$ and choose larger values of $d$, or keep the same $d$ and choose smaller $a$, we would eventually reach a point where there would be only a stable limit cycle with small $Y$.  If we were to instead choose smaller values of $d$, we would first move to a region where the small $Y$ equilibrium becomes stable and then to a region where only the large $Y$ equilibrium is present.


\section{Discussion and Conclusions}\label{sec:discussion}

In this study, we created and analyzed a relatively simple model that incorporates attitude change toward vaccination in response to the current disease prevalence.  Asymptotic approximation allowed us to obtain simple conditions for the existence and stability of disease-free and endemic disease equilibria.  We then chose some example functional forms for the rate constants in the attitude change processes to provide examples.

The closest prior work to ours is the collection of papers that use an information index to modify the vaccination rate \cite{bauch2005, buonomo2022, buonomo2013, donofrio2012, donofrio2011, donofrio2007, zuo2023}.  The two post-COVID papers in this group \cite{buonomo2022, zuo2023} lack a complete analysis of endemic disease equilibria.  The earlier papers contain a complete EDE analysis; in each case, a larger sensitivity of vaccination to the information index can lead to an unstable EDE with a stable limit cycle.  The primary differences between the information index models and our flexible attitude model are the mechanism by which changes in prevalence impact vaccination rates and the incorporation in our model of loss of immunity.

Even when the information index is based on recent prevalence history rather than just current prevalence, vaccination rates based on it respond much faster to prevalence changes than they do in our model.  This is because of the different assumptions about human behavior.  Prior to COVID-19, it was natural to assume that decisions about vaccination were based on a rational comparison of the payoffs of vaccination as compared to non-vaccination.  Rational choices happen quickly, on the time scale of the disease processes.  Following COVID-19, it is clear that vaccination decisions are based more on individual attitudes toward vaccination rather than rational comparison of alternative strategies.  Unlike rational assessment, attitude changes are slow because people tend to be resistant to changing their opinions in response to new facts.  In our model, vaccination rates respond to changes in prevalence on a slower time scale because the mechanism is indirect, with prevalence determining rates of attitude change rather than attitude itself and rate changes occurring on a slower time scale than that of prevalence changes.

Another implication of the difference between rational choice and attitude is that it makes sense to prescribe a mechanistic model for rational choice, whereas mechanistic models are less clear for human behavior.  While it would be natural to expect an increase in a pro-vaccination attitude as disease prevalence increases, we also considered the possibility of a non-monotone response, in which a sufficiently high prevalence eventually results in loss of confidence in the vaccine, thereby resulting in a decrease in pro-vaccination attitude with increasing prevalence.

For the monotone case, there is always a unique endemic disease equilibrium whenever the disease-free equilibrium is unstable.  This unique EDE is sometimes stable, but when the response of the rate change coefficient to disease prevalence is strong enough, the EDE becomes unstable and the resulting behavior converges to a limit cycle.

The non-monotone case always has at least one EDE, but it can have up to three in the example function family we used, and possibly more for an unrealistically complicated function family.  As seen in Figure \ref{fig:Bifplots2}, a variety of possibilities for stability exist, including either a stable EDE or a stable limit cycle when there is a unique EDE and an unstable middle EDE with a stable large EDE and either a stable small EDE or a stable limit cycle when there are three EDE's.

While the model is speculative in nature, it does suggest the possibility that attitude changes in response to prevalence can destabilize an endemic disease equilibrium, with multiple long-term outcomes in the case where the response is non-monotone.  It would perhaps be difficult to identify this type of instability in real data because it would have seasonal and stochastic variation superimposed on it, but it would be worth exploring.

In searching for real examples of EDE destabilization, it is helpful to know what disease features are likely to promote destabilization.  Most obviously, the sensitivity of the rate constants for attitude changes to differences in prevalence should be key, as evidenced by the greater tendency toward instability in the parameters $\Sigma$, $a$, and $d$ in the functions we used for these rate constants.  It is also worth noting the destabilizing influence of the rate constant for loss of immunity $h$.  While we might expect faster loss of immunity to be destabilizing, the plots of Figure \ref{fig:EDEplots1} suggests otherwise, as it shows that a larger value of the attitude variability parameter $\Sigma$ is needed for instability for a relatively large value of $h$ as compared to a disease where there is no loss of immunity.  But note also from this figure that the overall incidence of the disease is much higher for diseases with short-lived immunity than diseases with more permanent immunity, as of course would be expected.

The use of asymptotic methods to produce approximate analytical results has been studied in detail by Ledder \cite{ledder2023using}, with the conclusion that the actual difference in prediction between results for $\epsilon \to 0$ and small values like $\epsilon=0.0005$ are insignificant compared to differences caused by small changes in the estimates used for the model parameters.  Experiments using eigenvalues to determine stability for parameter sets not extremely close to the region boundaries confirm the same result here.

In an effort to focus on the process of attitude change motivated strictly by disease incidence, we made a number of simplifying assumptions: simple SIR disease class structure, perfect vaccination that prevents infection and (through boosters) loss of immunity, no influence of incidence history, and no drivers of attitude connected to interpersonal or media-based communication.  Of these assumptions, the assumption of perfect vaccination is the most likely to make a qualitative difference in outcomes, as imperfect vaccination would make prevalence less sensitive to vaccination rates and likely reduce the tendency toward instability.

While our study does not directly apply to COVID-19, it does suggest possibilities for how attitude changes will impact COVID-19 prevalence.  While the actual epidemiological structure for COVID-19 should include a latent (usually called `exposed') class and an asymptomatic class, these differences are unlikely to affect qualitative outcomes in the long run.  Gradual loss of immunity is particularly important for COVID-19 population dynamics, with a typical loss of immunity over 6 months corresponding to a much larger value of $h$ than we used in our generic disease examples.


\appendix

\section{Proposition \ref{DFE}}\label{appendix:DFEprop}

In this and the other appendices, we make use of the additional symbol definition
\beq
\label{params2}
\xi=\sigmabar+h\psibar=\hbar \sigmabar-h \omega,
\eeq
which is included in Table \ref{tb:grouped}.

With $Y = Z =0$,\footnote{Throughout this section, the functions $\omega$ and $\psi$ are evaluated at $Y=0$.} the remaining equilibrium values satisfy the equations
\begin{align}
\label{dfe1}
0 &= 1-S-vP+hX, \\
\label{dfe2}
0 &= \omega S-\zeta P, \\
\label{dfe3}
0 &= vP-R-hX, \\
\label{dfe4}
0 &= \psi R-\eta X.
\end{align}
Solution of the combination \eqref{dfe1} and \eqref{dfe3}, followed by \eqref{dfe2} and then \eqref{dfe4}, yields
\beq
R=1-S, \qquad P=\frac{\omega}{\zeta} S, \qquad X=\frac{\psi}{\eta}R.
\eeq

Substitution of these formulas into \eqref{dfe3}, and using the identity $\eta+h \psi = \xi$, ultimately yields the result
\beq
S = \frac{\zeta \xi}{\zeta \xi + v \omega \eta}.
\eeq
The vaccine-reduced reproduction number $\rv$ can be found from the scaled problem \eqref{syst2} using the next generation method (see \cite{brauer2019mathematical}, for example). We define the matrices $\mathcal{F}$ and $\mathcal{V}$ by
\[
\mathcal{F}=\begin{bmatrix}
\Gamma \rr SY \\
\Gamma \rr PY
\end{bmatrix},\quad \mathcal{V}=\begin{bmatrix}
    \Gamma Y \\
    \Gamma Z
\end{bmatrix},
\]
where $\Gamma=\epsilon^{-1}$: $\mathcal{F}$ and $\mathcal{V}$ represent the inflow and outflow of individuals to and from the disease compartments, respectively. We then compute
\[
F = \begin{bmatrix}
    \partial_Y\mathcal{F}_1 & \partial_Z\mathcal{F}_1\\
    \partial_Y\mathcal{F}_2 & \partial_Z\mathcal{F}_2
\end{bmatrix} = \begin{bmatrix}
    \Gamma \rr S & 0\\
    \Gamma \rr P & 0
\end{bmatrix}
\]
and
\[
V = \begin{bmatrix}
    \partial_Y\mathcal{V}_1 & \partial_Z\mathcal{V}_1\\
    \partial_Y\mathcal{V}_2 & \partial_Z\mathcal{V}_2
\end{bmatrix} = \begin{bmatrix}
    \Gamma & 0\\
    0 & \Gamma
\end{bmatrix},
\]
where $P$ and $S$ are evaluated at the disease-free equilibrium. The next generation matrix $K_L$ is given by
\[
K_L = FV^{-1} = \begin{bmatrix}
    \Gamma \rr S & 0\\
    \Gamma \rr P & 0
\end{bmatrix}
\begin{bmatrix}
    \epsilon & 0\\
    0 & \epsilon
\end{bmatrix} = \begin{bmatrix}
    \rr S & 0\\
    \rr P & 0
\end{bmatrix}.
\]
The vaccine-reduced reproduction number $\rv$ is given by the largest positive eigenvalue of $K_L$; hence,
\[ \rv = \rr S=
\frac{\zeta \xi}{\zeta \xi + v \omega \eta} \rr. \]

To determine the stability of the DFE, we begin with the Jacobian:
\beq
\label{JacDFE}
J_{DFE} =
\left( \begin{array}{cr|ccrr}
\Gamma (\rv -1) & 0 & 0 & 0 & 0 & 0 \\
\Gamma p & -\Gamma & 0 & 0 & 0 & 0 \\
\hline
-s\;\;\; & 0 & -1\;\;\; & -v\;\;\; & 0 & h \\
\;A & 0 & \omega & -\zeta\;\;\; & 0 & 0 \\
1 & 0 & 0 & v & -1 & -h \\
\;B & -1 & 0 & 0 &\psi & -\eta
\end{array} \right),
\eeq
where all state variables are evaluated at the equilibrium and all functions of $Y$ at $Y=0$, lower case letters are defined in \eqref{eqvars}, and
\beq
A = S \omega'-P \Sigma'-p, \qquad B= 1+R \psi'-X \Sigma'.
\eeq
The block structure of the matrix allows us to immediately identify two of the eigenvalues.  The remaining eigenvalues come from the characteristic polynomial for the lower right block, which factors as
\beq
P_{3456}(\lambda) = (\lambda +1)f(\lambda),
\eeq
where
\[ f(\lambda)=\lambda ^3+c_1 \lambda ^2 +c_2\lambda +c_3, \] with
\begin{align*}
&c_1 = \eta +\zetabar, \\
&c_2=\zeta +\eta \zetabar +v\omega+h\psi,\\
&c_3 = \eta \zeta +v \omega \eta+h\psi \zeta.
\end{align*}
From
\[ c_1>\eta+\zeta, \qquad c_2>(\zeta+v\omega)+h\psi, \]
we have
\[ c_1c_2>\eta(\zeta+v\omega)+h\psi\zeta=c_3. \]
Hence, the eigenvalues coming from the roots of $f$ all have negative real part, and stability is determined entirely by $\rv<1$.
Replacing $\xi$ by $\hbar \sigmabar-h \omega$ and rearranging yields a minimum requirement for $\omega$ given all other parameter values:
\[ \omega(0) > \omega_{cr} \equiv
\frac{ r\hbar \sigmabar \zeta}{rh \zeta+v \eta} =
\frac{ r\hbar \sigmabar \zeta}{\sigmabar(v+rh)+\rr hv}, \]
where $r=\rr-1$.

\section{Proposition \ref{EDEomegastar}}\label{appendix:omegaprop}

\subsection{Existence Condition for the Endemic Disease Equilibrium}

In the asymptotic limit $\epsilon \to 0$, the system \eqref{syst2} yields an equilibrium system in which the first two equations and the fourth reduce to
\[ s=1, \qquad z=yp, \qquad \chi p=\omega, \]
where
\beq
\chi=\zeta+y,
\eeq
after which the remaining equations are
\begin{align}
\label{ede1}
\rr-1-vp+hx-y &= 0, \\
\label{ede2}
y+vp-\rr R-hx=0\\
\label{ede3}
y-yp+\psi \rr R-\eta x &= 0.
\end{align}

\n Then \eqref{ede1} and \eqref{ede2} combine to yield
\[ \rr R=\rr-1=r \]
and \eqref{ede2} and \eqref{ede3} combine (replacing $\psi$ by $\Sigma-\omega$ and using $y=\chi-\zeta$) to yield
\[ \sigmabar x=\sigmabar r+\sigmabar p-\rr \chi p. \]
Substituting the last result into \eqref{ede2} yields an equation that can be solved for $p$ in terms of $y$ and $\Sigma$, along with $\rr$, $v$ and $h$:
\[ p = \frac{\sigmabar (r \hbar -y)}{\sigmabar (v+rh)+\rr h(v+y)}=\frac{\sigmabar (r \hbar -y)}{\rho+\rr hy}. \]
Using $\omega=\chi p$, this yields the equilibrium relation
\beq
\label{omegastar}
\omega(Y)=\omega^*(\rr Y;\Sigma, \rr, v, h), \qquad \omega^*(y)=\frac{\sigmabar (r \hbar -y)(\zeta+y)}{\rho+\rr hy}.
\eeq

\subsection{Properties of the Function $\omega^*$}

The properties $\omega^*(0)=\omega_{cr}$ and $\omega^*(r\hbar)=0$ follow immediately from the function definition.  For the derivative condition, we begin by writing the function as
\[ \omega^* = \sigmabar \;\frac{r\hbar \zeta-(\zeta-r\hbar) y-y^2}{\rho+\rr hy}, \]

Differentiating yields
\[ {\omega^*}'=-\sigmabar \; \frac{(\zeta -r\hbar) \rho+\rr h r\hbar \zeta+2\rho y+\rr hy^2}{(\rho+\rr hy)^2} =-\sigmabar \; \frac{\kappa+2\rho y+\rr hy^2}{(\rho+\rr hy)^2}. \]
Thus, $\kappa>0$ is sufficient for monotonicity.
Note that
\[
\zeta \rho = \zeta \left[ \left( \sigmabar+\rr h \right) v+\sigmabar rh \right] \ge v \left[ \left( \sigmabar+\rr h \right) \zeta+\sigmabar rh \right];
\]
hence,
\[ \kappa \ge v \left[ \left( \sigmabar+\rr h \right) \zeta+\sigmabar rh \right]-\sigmabar r\hbar v = v \left[ \left( \sigmabar+\rr h \right) \zeta-\sigmabar r \right]. \]
Replacing $r$ by $\rr-1$ yields
\[ \kappa \ge v \left[ \left( h \zeta-\sigmabar \right) \rr +\sigmabar \bar{\zeta} \right]. \]
We therefore have the result that
\[ (\sigmabar-h\zeta)\rr< \sigmabar \bar{\zeta} \]
is sufficient for monotonicity, albeit not necessary.
This requirement is satisfied for any reasonable set of parameters.  The most critical case has $h=0$ and $\Sigma \to 0$, in which case the condition reduces to $\rr<v+2$, which is almost certainly true.  Alternatively, $h\zeta>\sigmabar$ is clearly sufficient, and this reduces to $h>\sigmabar/(v+\sigmabar)$, with the latter likely small and certainly less than 1.

\section{Proposition \ref{EDEstability}}\label{appendix:EDEstab}

The Jacobian matrix at the endemic disease equilibrium is
\beq\label{JacEDE}
J_{EDE} =
\left( \begin{array}{crccrr}
0 & 0 & y\Gamma & 0 & 0 & 0 \\
p\Gamma & -\Gamma & 0 & y\Gamma & 0 & 0 \\
-1\;\;\; & 0 & -\ybar\;\;\; & -v\;\;\; & 0 & h \\
A & 0 & \omega & -u\;\;\; & 0 & 0 \\
1 & 0 & 0 & v & -1 & -h \\
B & -1 & 0 & 0 &\psi & -\eta
\end{array} \right),
\eeq
where
\beq
\label{eqAB}
\Gamma = \epsilon^{-1}, \qquad A = S \omega'-P \Sigma'-p, \qquad B = 1+R \psi'-X \Sigma'
\eeq
have been introduced to simplify the notation.

The usual way to compute the characteristic polynomial for a matrix $J$ is to find the determinant of $\lambda I-J$; however, it is more efficient to use the characteristic polynomial theorem \cite{ledder2023using}:
\begin{theorem}
For an $n \times n$ matrix $J$, let $S$ be the set of all nonempty subsets of the integers $1, 2, \ldots, n$.  For each possible $K \in S$, let $J_K$ be the determinant of the submatrix of $J$ that contains the entries in the rows and columns indicated by the index set $K$.  Then the characteristic polynomial of $J$ is
\beq
P(\lambda)=\lambda^n + c_1 \lambda^{n-1} + c_2 \lambda^{n-2} + \cdots + c_{n-1} \lambda + c_n,
\eeq
where
\beq
c_m=(-1)^m \sum_{|K|=m} J_K, \qquad c_n=(-1)^n |J|.
\eeq
\end{theorem}

\n While this theorem technically requires the computation of all possible subdeterminants of J, some can be ignored when the goal is merely to find leading order approximations of the coefficients.

For \eqref{JacEDE}, we can quickly get leading order results
\[ c_1 \sim \Gamma, \qquad c_2 \sim (\chi+\eta+2\ybar)\Gamma, \qquad c_3 \sim -J_{123}=y\Gamma^2, \]
where the last of these results is facilitated by the observation that the other $3 \times 3$ subdeterminants are only $O(\Gamma)$ or smaller.  To get $c_4$, we must calculate the three $4 \times 4$ subdeterminants of $O(\Gamma^2)$:
\[ J_{1234} \sim (vA+\chi)y \Gamma^2, \qquad J_{1235} \sim y\Gamma^2, \qquad J_{1236} \sim (\eta+hp-hB) y\Gamma^2; \]
hence,
\[ c_4 \sim [vA-hB+(\chi+hp+\etabar)] y \Gamma^2. \]

We next compute the full determinant to get
\[ c_6 = [(v\eta+hy)A-h\chi B+\chi(hp+\eta)] y \Gamma^2. \]
Finally, we need three of the $5 \times 5$ subdeterminants:
\[ J_{12345} \sim -(vA+\chi)y \Gamma^2, \qquad J_{12346} \sim -c_6, \qquad J_{12356} \sim -[(p-B)h+\eta] y\Gamma^2. \]
Combining these gives the convenient formula
\[ c_5 \sim c_6+c_4-c_3. \]
To leading order, we have the characteristic polynomial as
\beq
\label{charpoly}
P(\lambda)= \lambda ^6 + \Gamma \lambda ^5 +k_2 \Gamma \lambda ^4 + y\Gamma ^2\lambda ^3 +k_4\Gamma ^2 \lambda ^2 +k_5\Gamma ^2 \lambda +k_6\Gamma ^2,
\eeq
where
\beq
\label{k2k4}
k_2 = \bar{\chi}+\bar{\eta} + 2y, \qquad k_4 = [vA-hB+(\chi+hp+\bar{\eta})]y,
\eeq
\beq
\label{k5k6}
k_6 = [(v \eta+hy)A-h \chi B+\chi(hp+\eta)]y, \qquad k_5=k_6+k_4-y.
\eeq

Using only the leading order terms, we obtain the Routh array \cite{ledder2023using} as
\[
\begin{array}{rrrr}
1 &k_2 \Gamma &k_4 \Gamma^2 &k_6 \Gamma^2\\
\Gamma & y\Gamma^2 &k_5\Gamma ^2 &\\
q_1\Gamma &k_4\Gamma ^2 &k_6\Gamma ^2 &\\
\frac{q_2}{q_1} \Gamma ^2 & \frac{q_3}{q_1}\Gamma ^2 & &\\
k_4\Gamma ^2 &k_6\Gamma ^2 & &\\
\frac{q_4}{k_4} \Gamma ^2 & & &\\
k_6\Gamma ^2 & & &
\end{array},
\]
where
\begin{align*}
q_1 &= k_2-y=\bar{\chi}+\etabar+y>0, \\
q_2 &= yq_1-k_4=y [\ybar-hp-vA+hB], \\
q_3 &= k_5 q_1-k_6, \\
q_4 & = \frac{k_4 q_3-k_6 q_2}{q_1} = k_4k_5-yk_6 = (k_4+k_6)(k_4-y).
\end{align*}
Stability requires the quantities in the left column of the array to be positive.  Only three conditions need to be checked, as $q_1>0$ is immediate and $q_4>0$ is equivalent to $k_4>y$, which then guarantees $k_4>0$.  From $q_2>0$, $k_4>y$, and $k_6>0$, we obtain the conditions
\beq
\label{q2}
vA-hB < \ybar-hp,
\eeq
\beq
\label{k4}
hB-vA < \chi+hp+\eta,
\eeq
\beq
\label{k6}
h\chi B-(v\eta+hy)A < \chi(hp+\eta).
\eeq
Substituting in the formulas for $A$ and $B$ from \eqref{eqAB} yields the three criteria of Proposition \ref{EDEstability}.

Taking $\Sigma'=0$ reduces the stability criteria to
\begin{align}
\omega'<\rr \Upsilon_1(\rr Y), & \qquad \Upsilon_1 (y)=\frac{w_1}{v+rh}, \\
-\omega'<\rr \Upsilon_2(\rr Y), & \qquad \Upsilon_2 (y)=\frac{\sigmabar+w_2}{v+rh}, \\
-\omega'<\rr \Upsilon_3(\rr Y), & \qquad \Upsilon_3 (y)=\frac{\sigmabar w_2}{v\eta+rh\chi+hy},
\end{align}
where
\beq
w_1 = h(1-p)+vp+\ybar, \qquad w_2 = v(1-p)+hp+\ybar+\Sigma.
\eeq

Now suppose $\Upsilon_3>\Upsilon_2$, which would make the second criterion more limiting than the third.  Then
\begin{align*}
\sigmabar w(v+rh) &> (w+\sigmabar)(v\eta+rh\chi+hy) \\
&> w(v\eta+rh\chi) \\
&> w(\sigmabar v+\sigmabar rh);
\end{align*}
this contradiction establishes that the $\Upsilon_3$ criterion is always more limiting than the $\Upsilon_2$ criterion.

To obtain the result given in the last part of Proposition \ref{EDEstability}, we employ the identities
\[ v\eta+rh\chi+hy=v(\sigmabar+h)+ rh(\sigmabar+v+y)+hy=\rho+\rr hy \]
and
\[ w_2 = (\zeta+y)-(v-h)\frac{\sigmabar (r\hbar-y)}{\rho+\rr hy} = \cdots = \frac{\kappa+2\rho y+\rr hy^2}{\rho+\rr hy}. \]


\section*{Acknowledgments}
This material is based upon work that was done at the American Institute of Mathematics (AIM) with support from the National Science Foundation (NSF), and the authors thank AIM and NSF.

\bibliographystyle{abbrv}
\bibliography{bibliography}

\end{document}